\documentclass[11pt,reqno]{amsart}
\usepackage{graphicx}
\usepackage{amssymb,latexsym,epsfig,mathrsfs,graphpap,graphics,amsmath,amssymb}
\usepackage{amstext}
\usepackage{verbatim}
\usepackage{color}

\newtheorem{conj.}[thm]{Conjecture}

\theoremstyle{definition}

\theoremstyle{remark}

\numberwithin{equation}{section}

\begin{document}

\begin{flushleft}
  {\bf\Large {Embeddability, representability and universality \\ [1.5mm]  involving Banach spaces}}
\end{flushleft}

\parindent=0mm \vspace{.2in}

  {\bf{M. A. Sofi$^{\star}$}}

 \parindent=0mm \vspace{.1in}
{\small \it $^{\star}$Department of Mathematics, Central University of Kashmir, Srinagar-190015, Jammu and Kashmir, India. E-mail: $\text{aminsofi@gmail.com}$}

\parindent=0mm \vspace{.2in}
{\small {\bf Abstract:}  Given a category of objects, it is both useful and important to know if all the objects in the category may be realised as sub-objects - via morphisms in the given category- of a single object in that category enjoying some nice properties. In the category of separable Banach spaces with morphisms consisting of linear isometries, such an example of (a universal) object is provided by the well-known Banach Mazur theorem: the space $C[0,1]$ of continuous functions on the unit interval contains each separable Banach spaces as a closed subspace via a linear isometry.In other words, $C[0,1]$ is (isometrically) universal for the class of all separable Banach (even metric) spaces. Here the question also arises if, as opposed to realising (separable) Banach spaces as spaces of continuous functions on $[0, 1]$, it is possible to embed Banach spaces as subgroups of the group of linear isometries (resp. unitaries) on a nice Banach (resp. Hilbert) space. If such is the case, one says that the given Banach space is representable as a group of isometries (resp. unitaries).

\parindent=8mm \vspace{.1in}
On the other hand, the idea of embeddability involves the possibility of realising each object in a given class of objects as included inside another object of the same class enjoying some good properties which are not present in the initial object. A typical example of this phenomenon is provided by a well-known theorem of Pelczynski that each Banach space with the bounded approximation property may be embedded as a complemented subspace of a Banach space having a Schauder basis. Further, considering that a Banach space also comes equipped with weaker structures involving the underlying metric (Lipschitz), uniform and topological structures, it follows that besides the linear isomorphisms (isometries), one may also consider morphisms in this category consisting of maps which are Lipschitz, uniformly continuous or continuous. This motivates the consideration of situations where it becomes necessary to know if a Banach space (resp. a metric space) may be embedded in a nice Banach space as a metric, uniform or merely as a topological space.

\parindent=8mm \vspace{.1in}
The present paper deals with these and certain related issues and how they lead to some important questions in the nonlinear geometry of Banach spaces. In view of the restriction on the size of the paper as specified by the editorial board, no attempt will be made to provide detailed proofs of statements as included here, barring on a few occasions where proofs will be sketched in their barest minimum of details.

\parindent=0mm \vspace{.1in}
{\bf {Mathematics Subject Classification (2010):}} Primary 43A35, 22A25. Secondary 43A35.

\parindent=0mm \vspace{.1in}
{\bf { Key words and phrases:}} Hilbert space, Banach space, uniform embedding, representability.

\parindent=0mm \vspace{.2in}
{\bf Introduction}

\parindent=0mm \vspace{.1in}
The Banach Mazur theorem (see [9], Theorem 5.6) as referred to above which states that each separable Banach space embeds linearly and isometrically into $C[0,1]$ motivates the following natural questions:

\parindent=0mm \vspace{.1in}
{\bf a.} Is it possible to embed each (separable) Banach space as a subspace of a Hilbert space via an embedding which is Lipschitz or a uniform homeomorphism?

\parindent=0mm \vspace{.1in}
{\bf b.}	Is it possible to realise each Banach space (at least from a certain class of Banach spaces) into the space of unitaries of some Hilbert space (or, more generally, into the group of (linear) isometries of some nice Banach space)?

\parindent=0mm \vspace{.1in}
As we proceed, it will become clear that the above two questions are closely tied to each other.

\parindent=0mm \vspace{.1in}
{\bf 2.	Notation}

\parindent=0mm \vspace{.1in}
$\mathbb U$: Urysohn space

\parindent=0mm \vspace{.1in}
$\ell_2$: Space of square summable sequences

\parindent=0mm \vspace{.1in}
$\omega$:  Space of all numerical sequences

\parindent=0mm \vspace{.1in}
{\it Explanation:}

\parindent=0mm \vspace{.1in}

{\bf a.}	$\mathbb U$ is the unique complete, separable ultrahomogeneous metric space which contains isometric copies of all separable metric spaces.

\parindent=0mm \vspace{.1in}
Here by ``ultrahomogeneous", we mean that each isometry between finite subsets of $\mathbb U$  extends to an isometry on $\mathbb U$.

\parindent=0mm \vspace{.1in}
{\bf b.} $\ell_2=\Big\{ (x_n)\subset X: \sum_{n=1}^{\infty}\|x_n\|^2<\infty\Big\}$.

\parindent=0mm \vspace{.1in}
When equipped with the norm:
\begin{align*}
\pi_{2}\big((x_n)\big)=\left( \sum_{n=1}^{\infty}\|x_n\|^2\right)^{1/2},
\end{align*}
$\ell_2$ is a Hilbert space. (This space is also ultrhomogeneous in the above sense).

\parindent=0mm \vspace{.1in}
{\bf c.} $\omega=\big\{ (x_n): x_n\in \mathbb R\big\}=$ Countable product of the line.

\parindent=0mm \vspace{.1in}
Equipped with the product topology, $\omega$ is a nuclear Frechet space.

\parindent=0mm \vspace{.1in}
{\bf 3.	Morphisms in the category of Banach spaces}

\parindent=0mm \vspace{.1in}

Banach spaces may be identified/distinguished by means of one of the following five types of morphisms arising in Banach space theory:

\begin{enumerate}
  \item Isometric isomorphism
  \item Linear homeomorphism

  \item Lipschitz homeomorphism

  \item Uniform homeomorphism

  \item Topological homeomorphism
\end{enumerate}

\parindent=0mm \vspace{.1in}
In other words, one may identify Banach spaces X and Y according as there is a bijective map $T:X\to Y$ which satisfies (1), (2), (3), (4) or (5). Besides the structures determined by these morphisms, the one pertaining to the so called ‘coarse maps’ has been another important area of study in the theory of Banach spaces which has received considerable traction in recent years. However, we shall make no effort to discuss this particular aspect of the theory, and so have consciously avoided reference to the category with morphisms consisting of coarse mappings. Further, whereas (1) is the strongest of these morphisms, (5) is obviously the weakest. In fact, it is an old result of M. I. Kadec and Anderson (see [6], Chapter 6, Sections 9 and 10) that all infinite dimensional separable Banach spaces are mutually homeomorphic and that they are homeomorphic to $\omega$. Combining this fact with a theorem of V. Uspenskii [27], we have

\parindent=0mm \vspace{.1in}
{\bf Theorem 1.} {\it All the spaces introduced above are homeomorphic:}
\begin{align*}
\mathbb U\cong \ell_2 \cong \omega.
\end{align*}
The following consequences of Theorem 1 are worth noting which involves the homeomorphic  nature of the embeddings:

\parindent=0mm \vspace{.1in}
{\bf Corollary 2(a):} {\it Every separable metric space is homeomorphicic to a subspace of a Hilbert space.}

\parindent=0mm \vspace{.1in}
{\bf (b).} {\it Every separable metric space is homeomorphic to a subspace of a nuclear Frechet space.}

Regarding the other extreme in the above list involving isometric embeddings, we recall the following counterpart of the Banach Mazur theorem stated above:

\parindent=0mm \vspace{.1in}
{\bf Theorem 3 (see [9],Theorem 5.4]).} {\it Every separable metric space embeds isometrically into $\ell_{\infty}$, the Banach space of all bounded sequences (of scalars). For (real) normed spaces, the embedding may be chosen to be a linear isometry}.

\parindent=0mm \vspace{.1in}
{\bf Remark:} For finite metric spaces consisting of $n$ points, one can choose $\ell_{\infty}^n$ in place of $\ell_{\infty}$. With a little bit of effort, one can even use $\ell_{\infty}^{n-2}$ for $\ell_{\infty}^n$ as long as $n\ge 4$.

\parindent=0mm \vspace{.1in}
The latter part of the above theorem which is proved independently of the first part may be compared with the following striking theorems:

\parindent=0mm \vspace{.1in}
{\bf Theorem 4 (Mazur-Ulam).} {\it A bijective isometry $T:X\to Y$ acting between real normed spaces $X$ and $Y$ is already a linear isometry, provided $T(0)=0$. In particular, real normed spaces which are isometric as metric spaces (via an isometry fixing the origin) are actually isometric as normed spaces.}

\parindent=8mm \vspace{.1in}
There are several proofs of this statement having been devised ever since it was proved for the first time by these authors in 1932. However, the following short and crisp proof due to Nica [19] is simply irresistible which we include here for the sake of completeness. Indeed, in view of the continuity of $T$, it suffices to show that
\begin{align*}
T\left(\dfrac{x+y}{2}\right)=\dfrac{T(x)+T(y)}{2},\quad \text{for all}~x,y\in X.
\end{align*}
To this end, fix $x,y\in $ and let
\begin{align*}
\text{def}(T)=\left\|T\left(\dfrac{x+y}{2}\right)-\dfrac{T(x)+T(y)}{2}\right\|.
\end{align*}
We observe that
\begin{align*}
\text{def}(T)\le \dfrac{1}{2}\left\|T\left(\dfrac{x+y}{2}\right)-T(x)\right\|+\dfrac{1}{2}\left\|T\left(\dfrac{x+y}{2}\right)-T(y)\right\|=\left\| \dfrac{x-y}{2}\right\|.\tag{*}
\end{align*}

\parindent=0mm \vspace{.1in}
Clearly, $\text{def}(T)$ depends upon $T,x$ and $y$ and is bounded by a quantity which depends only upon $x$ and $y$ and is independent of $T$ as well as of $Y$. Now let the isometry $S$ on $X$ be given by $S(z)=T(x)+T(y)-z$ and consider the isometry $U=T^{-1} ST$ on $X$. We observe that $S(x)=y,S(y)=x$  and a simple calculation yields that def$⁡(U)=2\text{def}⁡(T)$. Finally, assuming
that $\text{def}(T)\ne 0$ and using iteration, we get an isometry on $X$ for which
$\text{def}(T)$ is arbitrarily large, contradicting the bound in $(*)$.

\parindent=0mm \vspace{.1in}
In the absence of surjectivity of the map in Theorem 4, we have

\parindent=0mm \vspace{.1in}
{\bf Theorem 5 (Godefroy-Kalton [12]).} {\it Let $X$ and $Y$ be separable Banach spaces and $T:X\to Y$ be an isometry. Then, $X$ is linearly isometric to a closed subspace of $Y$.}

\parindent=0mm \vspace{.1in}
A small digression: From Theorem 3, it follows, in particular, that each finite dimensional Banach space can be isometrically embedded into $C[0,1]$ and that obviously, no finite dimensional space can include (isometrically) all the finite dimensional Banach spaces as its subspaces. This motivates the following natural question involving universal spaces for the class of all $2$-dimensional spaces.

\parindent=0mm \vspace{.1in}
{\bf Question:} Does there exist a finite dimensional Banach space which is (isometrically) universal for all $2$-dimensional spaces?

\parindent=0mm \vspace{.1in}
The surprising answer to the previous question is no, as was shown by Bessaga [5]. Besides of course the space $C[0, 1]$, there exist ‘smaller’ spaces with this universal property. One such example is provided by $L_{1} [0,1]$. (For an elegant proof, see Yost [30]).

\parindent=0mm \vspace{.1in}
{\bf Remark 6:} The question arises: what happens if in place of (linear) isometries in the previous examples, one were to use linear mappings  $T:X\to Y$ which are $c$-isometries, i.e. such that $c^{-1} \|x-y\|\le \big\|Tx-Ty\big\|\le c\|x-y\|$ for $x,y\in X$ where $c>1$. Again surprisingly, the answer to this question now turns out to be in the affirmative under this weaker condition,  i.e., it is possible to choose a finite dimensional space which  contains  $c$-isomorphic copies of all  $2$-dimensional spaces! Indeed, using the compactness of the set of $2$-dimensional spaces Dim$(2)$ in the Banach Mazur distance, there exists a finite subset $S$ of  Dim$(2)$ such that any space in Dim$(2)$ is at a distance of at most $c$ from one of the spaces in $S$. It follows that the direct sum of the spaces in $S$ is a finite dimensional space which contains $c$-isomorphic copies of each $2$-dimensional space.

\parindent=0mm \vspace{.1in}
{\bf Using Lipschitz maps as morphisms}

\parindent=0mm \vspace{.1in}
Apart from $C[0,1]$  and $\ell_{\infty}$ which are universal for separable metric spaces in the sense of Theorem 2, it was shown by Aharoni [1] that the space $c_{0}$ (of null sequences in its sup-norm) provides yet another example of a (separable) universal Banach space, provided embedding is effected by means of a bi-Lipschitz mapping (See [4], Theorem 7.11 for a proof). Also, note that $c_{0}$ is not isometrically universal for separable metric spaces. In fact, it can be shown that $c_{0}$  does not contain an isometric copy of the unit circle in the plane. As pointed out in the latter reference, it amounts to saying that ``the space $c_{0}$ which is `small' in the linear category is quite `large' in the Lipschitz category".

\parindent=0mm \vspace{.1in}
{\bf Uniform Classification:}

\parindent=0mm \vspace{.1in}
As opposed to the existence of a (surjective) isometry which, as noted above, is strong enough to capture the linear structure of a Banach space, the uniform classification of Banach spaces, which deals with the study of Banach spaces as a uniform space, presents considerable difficulties that are not encountered in the isometric classification.

\parindent=8mm \vspace{.1in}
The metric induced by the norm of a Banach space defines a uniform structure on the given Banach space. Though it is obviously weaker than the structure induced by the norm, it is already strong enough to determine the linear structure of a Banach space. The following result provides a sample of these results.

\parindent=0mm \vspace{.1in}

{\bf Theorem 7.} {\it A locally convex space uniformly homeomorphic to a normed space is normable.}

\parindent=0mm \vspace{.1in}

{\it Proof (Sketch)}. Let $f:X\to Y$ be a uniform homeomorphism where $X$ is a normed space and $Y$ is locally convex such that $f(0) = 0$. Let $B$ be the unit ball in $X$. Then the set $U = f(B)$ is an open neighbourhood of $0$ in $Y$. By Kolmogorov’s theorem, it suffices to show that $U$ is bounded in $Y$, i.e. it is absorbed by each $0$-neighbourhood $V$ in $Y$, which by local convexity, can be chosen to be convex. Since $f$ is uniformly continuous, there exists a $0$-neighbourhood $W$ in $Y$ such that $x-y\in W$ implies that $f(x)-f(y)\in V$. Now there exists an $n$ such that each $x\in B$ can be written as $x=\sum_{i=1}^{n}w_{i}$ with  $w_{i}\in W$ for all $i$. This yields that $f\left( \sum_{i=1}^{m+1}w_{i}\right)-f\left( \sum_{i=1}^{m}w_{i}\right)\in V$ for all $1\le m< n$ and so, $f(x)=V+V+V+
 \dots+V=nV$, by convexity of $V$. In other words, $U=f(B)\subset nV.$

\parindent=0mm \vspace{.1in}
{\bf Corollary 8:} A nuclear Frechet space can never be uniformly homeomorphic to a Banach space, unless both spaces are finite dimensional (and thus having the same dimension).

\parindent=0mm \vspace{.1in}
The previous corollary suggest the following

\parindent=0mm \vspace{.1in}
{\bf Conjecture:} There does not exist a nuclear Frechet space which is uniformly universal for all separable Banach (metric) spaces.

\parindent=0mm \vspace{.1in}

The Hilbert space analogue of the above conjecture was already settled by the following well known theorem of Enflo.

\parindent=0mm \vspace{.1in}
{\bf Theorem 9([8]).} {\it A Banach space uniformly homeomorphic to a subspace of a Hilbert space is already a Hilbert space.}

\parindent=0mm \vspace{.1in}
The question involving the choice of a subset as opposed to a subspace of a Hilbert space as in Theorem 9 above  is interesting in its own right, leading to the notion of unitary representability as discussed in the next section.

\parindent=0mm \vspace{.2in}
{\bf 4.	Representing Banach spaces as groups of isometries}

\parindent=0mm \vspace{.1in}
{\bf Definition 10.} {\it A topological group $G$ is said to be unitarily representable if there exists a Hilbert space and a faithful continuous representation $\tau:G \to Iso(H)_s$. }

\parindent=0mm \vspace{.1in}
$G$ is said to be reflexively representable if in the above definition, $H$ is replaced by a reflexive Banach space $X$.

\parindent=0mm \vspace{.1in}
Equivalently, $G$ is unitarily representable if it is embeddable as a topological subgroup of $Iso(H)_s$. (Equivalently of $Iso(H)_w$). Similarly for reflexive representability.

\parindent=0mm \vspace{.1in}
{\bf Examples 11(a)}(See [15], Chapter 3). Every locally compact group G (and hence every finite dimensional Banach space) is unitarily representable. (Consider the right translations in $L_2 (G)$ given by $h\to g_h$ where $g_h (k)=kh,\,g,h\in G)$.

\parindent=0mm \vspace{.1in}
{\bf (b).} A Hilbert space is unitarily representable. This follows as a consequence of an important theorem of Uspenskij [28] asserting that the group of all rigid motions (including translations) of $\ell_2$ embeds isometrically as a subgroup of $Iso(\ell_2)_s$.

\parindent=0mm \vspace{.1in}
{\bf(c).}([3]) Every (real) nuclear space $X$ admits a unitary representation, which is faithful if $X$ admits a continuous norm.

\parindent=0mm \vspace{.1in}
{\bf (d).}([11]) Every Frechet Schwartz space is reflexively representable. Also, there exist Schwartz spaces which are not unitarily representable.

\parindent=0mm \vspace{.1in}
{\bf Remark 12}. There exist topological groups which are not reflexively representable. Megrelishvili [17] showed that $H_+ [0,1]$, the group of orientation preserving homeomorphisms on $[0,1]$ in its compact open topology is not reflexively representable. Consequently, $Iso(\mathbb U)$, the group of isometries on the Urysohn space is not reflexively representable (it contains $ H_+ [0,1]$ as a subgroup).

\parindent=0mm \vspace{.1in}
Amongst the classical Banach spaces, we have

\parindent=0mm \vspace{.1in}
{\it Example 1}: $L_p$--spaces for $1<p<\infty $ are reflexively representable.

\parindent=0mm \vspace{.1in}
{\it Example 2}: $L_p$--spaces for $1<p\le 2$ are unitarily representable. More generally, it follows from Uspenskii [28] that for the indicated range, the isometry group of $L_p$-spaces is unitarily representable.

\parindent=0mm \vspace{.1in}
{\it Example 3}: For $2<p<\infty $,$L_p$--spaces are not unitarily representable.

\parindent=0mm \vspace{.1in}
However, for arbitrary topological groups, we have

\parindent=0mm \vspace{.1in}
{\bf (c).} Every topological group $G$ admits a representation on a suitable Banach space which may be chosen to be $C_b^u (G)$, the space of bounded (right) uniformly continuous functions on $G$ with the sup-norm. (This is an old result of Teleman).

\parindent=0mm \vspace{.1in}
{\bf (d).} A unitarily representable group is always reflexively representable. However, converse is not true (See below).

\parindent=0mm \vspace{.1in}
Back to Enflo's theorem:

\parindent=0mm \vspace{.1in}
A nonlinear version of Enflo's theorem involving embedding into a subset instead of a subspace of a Hilbert space, is provided by the following theorem of Megrelishvili:

\parindent=0mm \vspace{.1in}
{\bf Theorem 13 ([17]).} {\it A Banach space is unitarily representable if and only if it embeds uniformly homeomorphic into a Hilbert space.}

\parindent=0mm \vspace{.1in}
{\it Proof (Sketch)}. The first part follows from the observation that $Iso(l_2)_s$ embeds uniformly into the Hilbert space direct sum
$\big(\sum_{n=1}^{\infty}(\ell_2)_n\big)_{\ell_2}$ via the map $g\to (gx_n)$ where the sequence $(x_n)$ can be chosen in $\ell_2$ such that $ \left \|x_n \right \| =2^{-n}$ and such that the maps :$ g\to gx_n$ generate the (left) uniformity of $Iso(\ell_2)_s$. Converse follows by using the fact that unitary representability of an Abelian group $G$ is equivalent to the condition that positive definite functions on $G$ separate the identity $e$ and closed subsets of $G$ not containing $e$. Combining it with Theorem 3.1 of [2] to the effect that uniform embeddability into a Hilbert space yields the latter condition completes the argument.

\parindent=0mm \vspace{.1in}
Back to Corollary 2:

\parindent=0mm \vspace{.1in}
$\ell_2$ is homeomorphically universal for all separable metric spaces. This motivates the following question.

\parindent=0mm \vspace{.1in}
{\bf Q1}: Is it true that  $\ell_2$ is also uniformly universal for the class of all separable metric spaces?
Before we answer this question, let us note that $l_2$ is not isometrically universal even for finite metric spaces. A simple example is provided by the (bipartite) graph $K_1, _3$ with respect to its graph distance. In fact, a finite simple connected graph can be isometrically embedded into $\ell_2$ if and only if it either a complete graph $K_n$ or a path $P_n$ for some $n$. However, finite metric spaces do admit bi-Lipschitz embeddings into a (finite dimensional) Hilbert space. Indeed, according to a celebrated theorem of Bourgain, an $n$-point metric space embeds into $\ell_2^n$ with distortion bounded by $c\log (n+1)$ where $c$ is a universal constant (See [20], Chapter 3).

\parindent=0mm \vspace{.1in}
Regarding Lipschitz universality of $\ell_2$, we have

\parindent=0mm \vspace{.1in}
{\bf (a).} ([20], Theorem 9.44): {\it The infinite binary tree is not Lipschitz embeddable into $l_2$.}

\parindent=0mm \vspace{.1in}
{\bf (b).} ([13], Section 12): {\it Heinsenberg group $H_3$ is not Lipschitz embeddable into $\ell_2$. 	However, as a locally compact group, $H_3$ is uniformly embeddable into $\ell_2$.}

\parindent=0mm \vspace{.1in}
Back to Question 1:

\parindent=0mm \vspace{.1in}
It turns out that the answer to $Q1$ is {\bf NO!}

\parindent=0mm \vspace{.1in}
In fact, Enflo (in answer to a question raised by Y. Smirnov) gave an example of a countable metric space which doesn't embed uniformly into $\ell_2$. However, we can still ask:

\parindent=0mm \vspace{.1in}
{\bf Q2}: Does there exist a (separable) reflexive Banach space which is uniformly universal for all (separable) complete metric spaces?

\parindent=0mm \vspace{.1in}
{\bf Remark 14(a)}. The non-existence of a separable reflexive Banach space which is isometrically universal for all separable reflexive Banach space is an old result of Szlenk [25].

\parindent=0mm \vspace{.1in}
The negative answer to $Q.2$ is contained in the following deep theorem of Kalton:

\parindent=0mm \vspace{.1in}
{\bf Theorem 15 ([14]).} {\it $c_0$ does not uniformly embed into a reflexive Banach space.}

\parindent=0mm \vspace{.1in}
{\bf Note (a)}: The Lipschitz analogue of the above assertion was established by Mankiewicz (1982).

\parindent=0mm \vspace{.1in}
Consequences:

\parindent=0mm \vspace{.1in}
{\bf (a).} $c_0$ is not reflexively representable. (Reflexive representability implies uniform embedding into a reflexive Banach space).

\parindent=0mm \vspace{.1in}
{\bf (b).} $C[0,1]$ is not reflexively representable. $C[0,1]$ does not even uniformly embed into a reflexive Banach space. (By Banach-Mazur theorem combined with $(a)$ above).

\parindent=0mm \vspace{.1in}
This motivates the following problem:

\parindent=0mm \vspace{.1in}
{\bf Problem 1.} Describe the class of groups/Banach spaces $X$ for which uniformly embeddability of $X$ into a reflexive Banach spaces is equivalent to $X$ being reflexively representable.

\parindent=0mm \vspace{.1in}
The difficulty involved in a solution to the  previous problem can be gauged from the following theorem of Yaacovet  al:

\parindent=0mm \vspace{.1in}
{\bf Theorem 16 ([29]).} {\it There exist a reflexive Banach space which is not reflexively representable.}

\parindent=0mm \vspace{.1in}
Such examples as guaranteed by Theorem 16 have to be located among no-classical Banach spaces as was shown to be the case for Tsirelson’s space and other Tsirelson-like spaces by these authors.

\parindent=0mm \vspace{.1in}
{\bf Remark 17.} From an important theorem of Aharoni, Maurey and Mityagin [2] combined with Theorem 13, it follows that a Banach space $X$ is unitarily representable if and only if it embeds isomorphically into $L_0 (\mu)$, the topological vector space of all complex-valued measurable functions on a measure space for some finite measure $\mu$, equipped with the topology of convergence in measure. Combining this with the fact that every Banach space isomorphic to a subspace of $L_0 (\mu)$ has cotype 2, it follows that a Banach space which is unitarily representable has cotype 2. Whether the converse is true is answered in the negative by combining a host of highly nontrivial results in Banach space theory proved over a period of time.

\parindent=0mm \vspace{.1in}
{\bf Example 18.}  There exist Banach spaces of cotype 2 which are not unitarily representable.

\parindent=0mm \vspace{.1in}
According to a famous theorem of N. T. Jaegermann [26], the Schatten class $S_p$ for $1<p<2$ has cotype 2 and type $p$. However, it does not embed into $L_0$ $(\mu)$. For otherwise, as a Banach subspace of $L_0 (\mu)$ having type $p,S_p$ can be isomorphically embedded into  $L_r$ for some $1<r<p$. But that would contradict an old theorem of Pisier [21] to the effect that $S_p$ with the range of $p$ as indicated above, cannot be isomorphically embedded into a Banach lattice with nontrivial cotype. Thus, $S_p$ is not unitarily representable.

\parindent=8mm \vspace{.1in}
This example may be compared with the case of the classical Lebesgue spaces $L_p$ which are unitarily representable for indicated range of $p$, as noted in Example 3 above.

\parindent=8mm \vspace{.1in}
The above discussion motivates the search for a class of cotype 2 Banach spaces which are unitarily representable. In this connection, it is useful to recall the definition of a class of Banach spaces having the so called Sazonov property which, as it turns out, can be embedded as subspaces of $L_0 (\mu)$, and so have the cotype 2 property. It appears likely that these spaces are unitarily representable.

\parindent=0mm \vspace{.1in}
{\bf Definition 19.} A Banach space $X$ is said to have the Sazonov property if there exists a locally convex topology $\tau$ on $X$ (weaker than the norm topology) such that each positive definite function $p$ on $X$ is continuous with respect to $\tau$ if and only if it can be expressed  as the Fourier transform against a positive Borel measure $\mu$ on $X^*$:
\begin{align*}
p(x)=\int e^{-i\left\langle x,x^{\ast}\right\rangle}d_\mu, \quad x\in X.\tag{*}
\end{align*}

\parindent=0mm \vspace{.0in}
{\bf Remark 20.} The requirement of a weaker topology on $X$ is dictated by the fact that on each infinite dimensional Banach space, there are positive definite functions which cannot be written as a Fourier transform as given above in $(*)$. Thus, whereas finite dimensional Banach spaces have the Sazonov property by virtue of Bochner's theorem for positive definite functions on locally compact groups, the property holds for Hilbert spaces $H$ where the desired topology is taken to be the projective topology determined by the quadratic forms: $ x\to \left \langle Tx,x \right \rangle $, where $T$ ranges over the set of all trace class operators on $H$ which are symmetric and positive. (See [7], Chapter 6).

\parindent=0mm \vspace{.1in}
{\bf Problem 2.} Investigate the relationship, if any, between unitary representability and the Sazonov property.

\parindent=0mm \vspace{.1in}
We conclude with the following problem which remains open:

\parindent=0mm \vspace{.1in}
{\bf Problem 3.} Do there exist non-abelian groups which are uniformly embeddable into a Hilbert space but are not unitarily representable?

\parindent=0mm \vspace{.2in}
{\bf 5.	Epilogue}

\parindent=0mm \vspace{.1in}
In the context of the Banach Mazur theorem which has been the main source of questions as discussed in the previous sections, there is yet another aspect of it which pertains to the question whether it could be achieved that (separable) Banach spaces $X$ may be embedded into $C[0,1]$ by means of a linear isometry $T$ so that $T(X)$ consists of functions enjoying better regularity properties than those which are merely continuous. Surprisingly, it turns out that this cannot be done if $T(X)$ is chosen to sit inside the (sub)space of (everywhere) differentiable functions unless, of course, $X$ is finite dimensional. This was observed by $W$. Lusky way back in the seventies.

\parindent=8mm \vspace{.1in}
Let us denote by $CND[0,1]$ the set of all nowhere differentiable continuous functions and by $D[0,1]$, the space of all everywhere differentiable functions on $[0,1]$.

\parindent=0mm \vspace{.1in}

{\bf Theorem 20(W. Lusky).} {\it $D[0,1]$ is not closed in $C[0,1]$. More generally, no infinite dimensional subspace of $D[0,1]$ can be closed.}

\parindent=0mm \vspace{.1in}
{\bf  Corollary 21.} {\it No (separable) infinite dimensional Banach space can be linearly and isometrically embedded into $D[0,1]$. In other words, a complete subspace of D[0,1] is necessarily finite dimensional.}

\parindent=0mm \vspace{.1in}
The above result provides the first hint of the ‘smallness’ of the (sub)space $D[0,1]$ which, as noticed above, turns out to be too small to contain all  separable infinite dimensional Banach spaces as subspaces. On the other hand, the following theorem shows that the set $CND[0,1]$ which is known to be nonempty- thanks to Weierstrass (see [24])-is actually quite large, both in the topological as well as in the category sense.

\parindent=0mm \vspace{.1in}
{\bf Theorem 22.} {\it $CND[0,1]$ is dense and of second category inside $C[0,1]$.}

\parindent=0mm \vspace{.1in}
Continuing in the same vein, the following result of Fonf, Gurariy and Kadec shows that $CND[0,1]$ is indeed  large, even in the `algebraic' sense!

\parindent=0mm \vspace{.1in}
{\bf Theorem 23([10]).} {\it $CND[0,1]$ contains a closed infinite dimensional vector subspace. Indeed, such a space would necessarily be uncountable dimensional!}

\parindent=0mm \vspace{.1in}
These considerations naturally lead to the following question:

\parindent=0mm \vspace{.1in}
{\bf Q3.} Is it true that $CND[0,1]$ is large enough that together with the identically zero function, it contains each separable Banach space linearly and isometrically (as a subset of $C[0,1]?$)

\parindent=0mm \vspace{.1in}
In the light of the above results which provide reasonable evidence to support an affirmative answer to the previous question, the following striking theorem of L. R. Piazza completes the story!

\parindent=0mm \vspace{.1in}
{\bf Theorem 24([23]).} {\it Every separable Banach space can be linearly and isometrically embedded into $CND[0,1]$ (inside $C[0, 1]$).}

\parindent=0mm \vspace{.1in}
The next question is motivated by a possible Lipschitz analogue of Theorem 17.

\parindent=0mm \vspace{.1in}
{\bf Q4.} Is it true that a closed subspace of $C[0,1]$ consisting entirely of Lipschitz functions is finite dimensional?

\parindent=0mm \vspace{.1in}
An affirmative answer to this question was provided by  M. Raja in 2015 ([22]).

\parindent=0mm \vspace{.1in}
We conclude with the following natural problems which to the best of our knowledge remain open.

\parindent=0mm \vspace{.1in}
{\bf Problem 4.} Is it true that each finite dimensional Banach space can be isometrically embedded into $D[0,1]?$

\parindent=0mm \vspace{.1in}
{\bf Problem 5.} Is it true that each finite dimensional Banach space can be isometrically embedded into $Lip[0,1]?$

\parindent=0mm \vspace{.1in}
As opposed to the Banach space case where, as we have seen, it is impossible to isometrically embed infinite dimensional spaces into spaces of everywhere differentiable functions, the case of Frechet spaces which are nuclear is far more satisfactory. In fact, it turns out that each nuclear Frechet space admits a linear isomorphic embedding into the space $C^\infty (\mathbb R)$ of infinitely differentiable functions on the line equipped with the (nuclear) Frechet topology given by the semi-norms:
\begin{align*}
\big\| f\big\|_{n}=\sup_{x\in U_{n},k\le n}\left| f^{(k)}(x)\right|
\end{align*}
where $U_{1}\subset U_{2}\subset \dots U_{n}\dots $ is an increasing sequence of open sets whose union is equal to $\mathbb R$. This follows from a famous theorem of $T$. Komura and Y. Komura ([18], Theorem 29.8) combined with an isomorphism theorem of D. Vogt (see [31], Theorem 3.22). An important takeaway from this result combined with Examples 11(c) and a host of similar results surrounding nuclear spaces is that they behave a lot more like finite dimensional spaces than do even infinite dimensional Hilbert spaces which lack, say the Heine-Borel property, contrary to the case of nuclear spaces where this latter property holds!

\end{document}